\documentclass[preprint,12pt]{elsarticle}
\usepackage[margin=1in]{geometry}
\usepackage[utf8]{inputenc}
\usepackage{graphicx}
\usepackage{multirow}
\usepackage{dcolumn}
\usepackage{bm}
\usepackage{color}      

\usepackage{amsmath}
\usepackage{amssymb}
\begin{document}
\begin{frontmatter}

\title{An eigenvalue problem for self-similar patterns in Hele-Shaw flows}

\author[1,2,3]{Wang Xiao}
\author[1,2,3]{Lingyu Feng}
\author[4]{Kai Liu \corref{cor2}}
\author[1,2,3]{Meng Zhao \corref{cor1}}
\address[1]{School of Mathematics and Statistics, Huazhong University of Science and Technology, Wuhan 430074, China}
\address[2]{Center for Mathematical Sciences, Huazhong University of Science and Technology, Wuhan 430074, China}
\address[3]{Steklov-Wuhan Institute for Mathematical Exploration, Huazhong University of Science and Technology, Wuhan 430074, China}
\address[4]{College of Education for the Future, Beijing Normal University, Zhuhai 519087, China}

\cortext[cor2]{Corresponding author: liuk@bnu.edu.cn}
\cortext[cor1]{Corresponding author: mzhao9@hust.edu.cn}



\begin{abstract}
Hele-Shaw problems are prototypes to study the interface dynamics. Linear theory suggests the existence of self-similar patterns in a Hele-Shaw flow. That is, with a specific injection flux the interface shape remains unchanged while its size increases. In this paper, we explore the existence of self-similar patterns in the nonlinear regime and develop a rigorous nonlinear theory characterizing their fundamental features. Using a boundary integral formulation, we pose the  question of self-similarity as a generalized nonlinear eigenvalue problem, involving two nonlinear integral operators. The flux constant $C$ is the eigenvalue and the corresponding self-similar pattern $\mathbf{x}$ is the eigenvector. We develop a quasi-Newton method to solve the problem and show the existence of nonlinear shapes with $k$-fold dominated symmetries. The influence of initial guesses on the self-similar patterns is investigated. We are able to obtain a desired self-similar shape once the initial guess is properly chosen. Our results go beyond the predictions of linear theory and establish a bridge between the linear theory and simulations.
\end{abstract}


\end{frontmatter}

\section{Introduction}
Various forms of pattern formation phenomena, ranging from the growth of bacterial colonies to the formation of snowflakes, exhibit analogous underlying physical mechanisms and mathematical structures. Understanding the formation kinetics and interplay of system parameters offers insight into pattern formation and improve control in various physical, biological, and engineering systems. 

 The Hele-Shaw flow is a classic problem for studying  interface dynamics. It is defined as a viscous flow between two parallel plates separated by a small gap \cite{hele1898flow, lee2011viscous, savina2011dynamical}.  Saffman and Taylor \cite{saffman1958penetration} found that injecting less viscous fluid into a viscous fluid initiates interface instabilities, resulting in the formation of a fingering pattern.  On the other hand, the surface tension restrains perturbation growth as suggested by the linear stability theory \cite{saffman1958penetration}.
 If we decompose the interface as a series of Fourier modes, the nonlinear interactions among different modes play a significant role and lead to the typical dense-branching form by repeated tip-splitting \cite{chuoke1959instability,langer1980instabilities,McLean1981TheEO,park1984two,maher1985development,ben1986formation, langer1989dendrites, ben1990formation, cummins1993successive, praud2005fractal,  li2007rescaling}. 
Brener {\it et al.} \cite{brener1990selection} conducted a study on viscous fingers in the $\pi/2$ sector geometry, suggesting that scaling the injection rate with time as $t^{-1/3}$ may enable the production of single fingers (no tip-splitting) with finite surface tension. This specific scaling factor is the dominant factor in time that allows for nonlinear self-similar evolution \cite{brener1990selection, combescot1991selection}. Combescot and Ben Amar \cite{combescot1991selection}, as well as Ben Amar {\it et al.} \cite{ben1991self}, numerically identify self-similar divergent (and convergent) fingers with finite surface tension. Li {\it et al.} \cite{PhysRevLett.102.174501} considered a radial Hele-Shaw cell and found that the limiting self-similar shape is actually universal. The limiting shape depends only on the flux constant and is independent of the initial configuration. Recently, numerous researchers \cite{2006Self, PhysRevFluids.4.064002, PhysRevFluids.7.053903, ZhaoSISC, PhysRevLett.115.174501, PhysRevE.103.063105, ZHAO2016394, zhao_niroobakhsh_2021, reinaud2022self} have conducted research on the self-similarity problem in Hele-Shaw cells under various conditions. Note that self-similar interface problems have been studied in material sciences. Li {\it et al.} \cite{li2004nonlinear} first proposed a nonlinear self-similar theory in a crystal growth problem. Barue {\it et al.} \cite{Amlan12} investigated the nonlinear simulations of the self-similar growth and shrinkage of a precipitate in inhomogeneous elastic media. Recently, the asymptotic and exact self-similar evolution of dendrite precipitate has been investigated in \cite{Amlan22}.

 Linear stability theories and nonlinear simulation results show the existence of self-similar patterns in the Hele-Shaw flow. 
 This paper aims to establish a rigorous nonlinear theory that goes beyond the linear theory and simulations, enabling a comprehensive study of the self-similar patterns in the radial Hele-Shaw problem. Using a boundary integral formulation, we derive the governing equations, a generalized nonlinear eigenvalue problem \cite{gilboa2018nonlinear}, $\displaystyle \mathcal{M}[{\bf x}] + C \mathcal{G}[{\bf x}] = 0$. These equations depend nonlinearly on the eigenvector ${\bf x}$, which represents the self-similar shape. The parameter $C$ is the critical flux constant and acts like an eigenvalue. Here $\mathcal{M}$ and $\mathcal{G}$ are nonlinear integral operators that incorporate high-order derivatives and a logarithmic singularity. 
This type of nonliear eigenvalue problems
appear in different applications in data sciences \cite{gilboa2018nonlinear}, physics \cite{bungert2022gradient,cohen2018energy}, mathematics \cite{bao2004computing, cances2010numerical, WeizhuBao2013Kinetic, henning2023dependency, dusson2023overview}, and so on. A different type of the Nonlinear Eigenvalue Problems (NEP)~\cite{guttel2017nonlinear}, where the underlying operator depends nonlinearly on the eigenvalue, has also been investigated.  See \cite{kublanovskaya1970approach, neumaier1985residual, kressner2009block, jarlebring2018disguised, lancaster1961generalised, schreiber2008nonlinear} for details.  


In this manuscript, $\displaystyle\mathcal{M}[\beta \mathbf{x}]= \beta^{-2}\mathcal{M}[\mathbf{x}]$ and $\displaystyle\mathcal{G}[\beta \mathbf{x}]=\beta\mathcal{G}[\mathbf{x}]$ for all $\beta\in R^+$. Thus, an eigenvector ${\bf x}$ can be paired with a given eigenvalue $C$ by changing its magnitude. To deal with the singular integrals in $\mathcal{M}$ and $\mathcal{G}$, we use the alternating trapzoidal rule \cite{greenbaum1993laplace}. Consequently, we implement a spectrally accurate quasi-Newton method to solve the nonlinear eigenproblem with a pre-specified $C$. Once $\mathbf{x}$ is obtained, we determine the eigenvalue $\displaystyle C = - \frac{\mathcal{M}[{\bf x}] }{ \mathcal{G}[{\bf x}]}$ under the same size of $\mathbf{x}$.


We find that the nonlinear flux of k-fold dominant self-similar shape ($k\ge4$) is smaller than that predicted by linear theories. However, the nonlinear flux of the 3-fold dominant self-similar shape exceeds the predictions of the linear theory. In addition, we find that the calculated self-similar shapes are closely related to the initial guesses. Surprisingly, we can obtain self-similar shapes that are not included in the initial guesses. A given self-similar shape may be obtained once the initial guess is properly chosen. This theory surpasses the realm of linear theory and achieves results comparable to simulations. The nonlinear theory serves as a bridge between linear theories and simulations.

This paper is structured as follows. In Section 2, we provide an overview of the governing equations and review the linear stability analysis. In Section 3, we present the formulation of the nonlinear self-similar theory. In Section 4, we introduce the quasi-Newton scheme for self-similar shapes. In Section 5, we present the numerical results and compare them to the predictions of the linear theory. In the final section,  we give conclusions and discuss future work.

\section{Review of the Hele-Shaw problem}
\subsection{Governing equations}
\begin{figure}[h]
    \centering
    \includegraphics[width=0.55\linewidth]{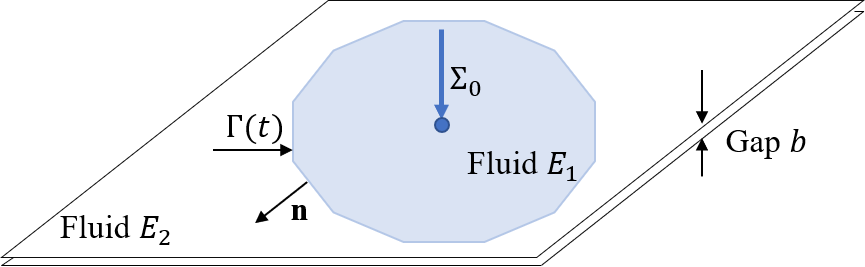}
    \caption{A schematic diagram for an air-oil interface system.}
    \label{fig1}
\end{figure}

We investigate a radial Hele-Shaw cell that features an air-oil interface system, as shown in Fig. \ref{fig1}. The system is comprised of one moving interface $\Gamma(t)$, which separates two fluid domains - a less viscous fluid domain $E_1$, and a more viscous fluid domain $E_2$. We assume fluids obey Darcy's Law,
\begin{equation}
    {\bf u}_i = -K_i \nabla P_i \quad \mathrm{in} \quad E_i \quad i=1,2, \label{velo}
\end{equation}
where ${\bf u}_i$ is the velocity of the fluid, $P_i$ is the corresponding pressure, and $K_i = b^2/(12\mu_i)$ is the mobility. $i=1, 2$ refers to the fluid inside and out of the interface, respectively. The parameter $b$ is the width of the gap between the two parallel plates of the Hele-Shaw cell and $\mu_i$ is the viscosity of fluids. We consider incompressible fluids and have $\nabla \cdot {\bf u}_i =0$. Thus, the pressure of the fluids satisfies
\begin{equation}
    \nabla^2 P_i = 0 \quad \mathrm{in} \quad E_i \quad i=1,2.
\end{equation}
Across the interface, the fluid normal velocity $V = {\bf u}_1 \cdot {\bf n} = {\bf u}_2 \cdot {\bf n}$ is continuous, where ${\bf n}$ is the outward normal. The pressure experiences a jump as dictated by the Young-Laplace condition,
\begin{equation}
    P_1 - P_2 = \tau \kappa \quad \mathrm{on} \quad \Gamma(t), \label{pressjump}
\end{equation}
where $\tau$ is the surface tension and $\kappa$ is the curvature of the interface. For simplicity, we consider an injection flux $J$ supplied at the origin, 
\begin{equation}
    J = \frac{1}{2\pi}\int_{\Sigma_0} \mathbf{u}_1\cdot \mathbf{n}ds, \label{flux}
\end{equation}
where $s$ is the arclength and $\Sigma_0$ is a small circle centered at the origin. The interface $\Gamma(t)$ evolves via
\begin{equation}
    {\bf n} \cdot \frac{d{\bf x}}{dt} = V \quad \mathrm{on} \quad \Gamma(t),
\end{equation}
where ${\bf x}$ is the position of the interface.

\subsection{Linear analysis}
Following \cite{saffman1958penetration, mullins1963morphological, paterson1981radial, buka1987stability}, we analyze the linear stability of the Hele-Shaw problem. Given an air-oil interface slightly perturbed by an azimuthal Fourier mode with wave number $k$, the interface evolve in the linear regime as
\begin{equation}
    r(\alpha, t) = R(t) + \epsilon\delta_k(t)\cos k\alpha,
\end{equation}
where $\epsilon\ll 1$, $\delta_k(t)$ is the amplitude of the perturbation,  $r(\alpha, t)$ and $R(t)$ are the radius of the perturbed and unperturbed air-oil interface, respectively. A classical linear stability analysis \cite{mullins1963morphological, paterson1981radial, buka1987stability} gives the growth rate of underlying circle $R(t)$ as
\begin{equation}
    R(t)\frac{dR}{dt} = J(t),
\end{equation}
and the shape factor $\displaystyle\frac{\delta_k}{R}(t)$ evolves as
\begin{equation}
    \left( \frac{\delta_k}{R} \right)^{-1}\frac{d}{dt}\left( \frac{\delta_k}{R} \right) = \frac{1}{R^2}(k-2)\left( J-\frac{C}{R} \right), \label{linear}
\end{equation}
where the constant $C=k(k^2-1)/(k-2)$ arises due to the stabilizing effects of surface tension. Equation (\ref{linear}) shows that the shape factor grows (decays) for $J>C/R$ ($J<C/R$). In particular, when $J=C/R$, the considered shape factor is time-independent, and the pattern is assumed to evolve self-similarly with a single mode perturbation under the constraints of linear theory. Therefore, using either a constant or increasing flux over time causes the growth of perturbations as $R$ increases. This leads to the development of ramified viscous fingering patterns~\cite{saffman1958penetration}.

\section{Nonlinear theory}
Within this section, we establish the complete set of nonlinear equations of a self-similarly evolving air-oil interface. Using potential theory \cite{kress1989linear,mikhlin2014integral}, the solutions of Eqs. (\ref{velo}) - (\ref{flux}) can be written in terms of boundary integrals. We represent the pressure on the interface $\Gamma$ through a single-layer potential and a double-layer potential, 
\begin{align}
    P_1({\bf x}) &=& \int_\Gamma \gamma_1({\bf x'}) G({\bf x}-{\bf x'})ds({\bf x'}) + \int_\Gamma \gamma_2({\bf x'}) \frac{\partial G({\bf x}-{\bf x'})}{\partial {\bf n}({\bf x'})}ds({\bf x'}) + \frac{1}{2}\gamma_2({\bf x}) - \frac{J}{K_1} \ln |{\bf x}|,\label{press1}\\
    P_2({\bf x}) &=& \int_\Gamma \gamma_1({\bf x'}) G({\bf x}-{\bf x'})ds({\bf x'}) + \int_\Gamma \gamma_2({\bf x'}) \frac{\partial G({\bf x}-{\bf x'})}{\partial {\bf
    n}({\bf x'})}ds({\bf x'}) - \frac{1}{2}\gamma_2({\bf x}) - \frac{J}{K_2} \ln |{\bf x}|.\label{press2}
\end{align}
Here, $G({\bf x}) = 1/2\pi\log|{\bf x}|$ is the Green's function. $\gamma_1({\bf x})$ and $\gamma_2({\bf x})$ represent the density on the interface for single- and double-layer potential, respectively. Across the interface, the pressure has a jump given by the Young-Laplace condition $P_1 - P_2 = \tau \kappa$, thus
\begin{equation}
    \gamma_2({\bf x}) = \tau \kappa + \frac{K_2 - K_1}{K_1 K_2} J\ln|{\bf x}|.\label{mu1}
\end{equation}

According to jump relations of potential theory \cite{kress1989linear} and Eqs. (\ref{press1}) - (\ref{press2}), we have the normal derivative of the pressure from interior and exterior region to the interface $\Gamma$,
\begin{align}
    \frac{\partial P_1({\bf x})}{\partial {\bf n}({\bf x})} &=& \int_\Gamma \gamma_1({\bf x'}) \frac{\partial G({\bf x}-{\bf x'})}{\partial {\bf n}({\bf x})} ds({\bf x'}) -\frac{1}{2}\gamma_1({\bf x}) + \int_\Gamma \gamma_2({\bf x}') \frac{\partial^2G({\bf x}-{\bf x}')}{\partial {\bf n}({\bf x})\partial {\bf n}({\bf x}')}ds({\bf x'})  - \frac{J}{K_1} \frac{{\bf x}\cdot{\bf n}}{|{\bf x}|^2}, \label{der1}\\
    \frac{\partial P_2({\bf x})}{\partial {\bf n}({\bf x})} &=& \int_\Gamma \gamma_1({\bf x'}) \frac{\partial G({\bf x}-{\bf x'})}{\partial {\bf n}({\bf x})} ds({\bf x'}) +\frac{1}{2}\gamma_1({\bf x}) + \int_\Gamma \gamma_2({\bf x}') \frac{\partial^2G({\bf x}-{\bf x}')}{\partial {\bf n}({\bf x})\partial {\bf n}({\bf x}')}ds({\bf x'})  - \frac{J}{K_2} \frac{{\bf x}\cdot{\bf n}}{|{\bf x}|^2}.\label{der2}
\end{align}
As the fluid normal velocity is continuous, ${\bf u}_1 \cdot {\bf n} = {\bf u}_2 \cdot {\bf n} = V({\bf x})$, and ${\bf u_i} = -K_i \nabla P_i$, we get $\partial P_1({\bf x})/\partial {\bf n}({\bf x}) = - V/{K_1}$ and $\partial P_2({\bf x})/\partial {\bf n}({\bf x}) = - V/{K_2}$. Combining Eqs. (\ref{der1}) and (\ref{der2}), we obtain
\begin{align}
    &\gamma_1({\bf x}) = \frac{K_2-K_1}{K_1K_2}\left(V - J\frac{{\bf x}\cdot{\bf n}}{|{\bf x}|^2} \right),\label{mu2}\\
    &\frac{1}{2}(K_1+ K_2)\gamma_1({\bf x}) + (K_2 - K_1) \int _\Gamma \gamma_1({\bf x}') \frac{\partial G({\bf x}-{\bf x'})}{\partial {\bf n}({\bf x})}ds({\bf x'}) = (K_1- K_2) \int_\Gamma \gamma_2({\bf x}') \frac{\partial^2G({\bf x}-{\bf x}')}{\partial {\bf n}({\bf x})\partial {\bf n}({\bf x}')}ds({\bf x}').\label{velocity}
\end{align}
Substituting Eqs. (\ref{mu1}) and (\ref{mu2}) into Eq. (\ref{velocity}), we get
\begin{equation}
    \begin{split}
    \frac{1}{2}(K_1+K_2)\left(V - J\frac{{\bf x}\cdot{\bf n}}{|{\bf x}|^2} \right) + (K_2-K_1)\int_\Gamma \left(V - J\frac{{\bf x}'\cdot{\bf n}}{|{\bf x}'|^2} \right) \frac{\partial G({\bf x}-{\bf x'})}{\partial {\bf n}({\bf x})}ds({\bf x'}) \\
    = -K_1 K_2 \int_\Gamma \tau \kappa({\bf x}') \frac{\partial^2G({\bf x}-{\bf x}')}{\partial {\bf n}({\bf x})\partial {\bf n}({\bf x}')}ds({\bf x}') + (K_1-K_2)J \int_\Gamma \ln|{\bf x}'| \frac{\partial^2G({\bf x}-{\bf x}')}{\partial {\bf n}({\bf x})\partial {\bf n}({\bf x}')}ds({\bf x}'). \label{solu}
    \end{split}
\end{equation}

Now we consider the self-similar shape. We can separate the time and space, 
\begin{equation}
    {\bf x} = R(t) \Tilde{\bf x}(s),
\end{equation}
where $\Tilde{\bf x}(s)$ is the self-similar shape and $R(t)$ is a scaling function only dependent on time. As a consequence, we know that $V({\bf x}) = (\Tilde{\bf x}\cdot {\bf n}) \dot{R}$, where the dot represents the derivative with respect to time. Also, we know that from the mass conservation, $R\dot{R} = \pi J/\Tilde{A}$, where $\Tilde{A}$ is the area enclosed by the self-similar shape $\Tilde{\bf x}$. Plugging into Eq. (\ref{solu}), we have
\begin{equation}
    \begin{split}
        \left[\frac{1}{2}(K_1+K_2) (1 - \frac{\Tilde{A}}{\pi |\Tilde{\bf x}|^2})\Tilde{\bf x}\cdot{\bf n} + (K_2-K_1) \int_{\tilde{\Gamma}} (1 - \frac{\Tilde{A}}{\pi |\Tilde{\bf x}'|^2}) \Tilde{\bf x}'\cdot{\bf n} \frac{\partial G(\Tilde{\bf x} - \Tilde{\bf x}')}{\partial {\bf n}(\Tilde{\bf x})}ds(\Tilde{\bf x}') \right] \dot{R}R^2 \\
        - (K_1-K_2)\frac{\Tilde{A}}{\pi} \int_{\tilde{\Gamma}} \ln{|\Tilde{\bf x}'|} \frac{\partial^2G(\Tilde{\bf x}-\Tilde{\bf x}')}{\partial {\bf n}(\Tilde{\bf x})\partial {\bf n}(\Tilde{\bf x}')}ds(\Tilde{\bf x}') \dot{R}R^2 = -K_1K_2 \int_{\tilde{\Gamma}} \tau\Tilde{\kappa}({\bf x}') \frac{\partial^2G(\Tilde{\bf x}-\Tilde{\bf x}')}{\partial {\bf n}(\Tilde{\bf x})\partial {\bf n}(\Tilde{\bf x}')}ds(\Tilde{\bf x}'). \label{equaG}
    \end{split}
\end{equation}

We define the following integral operators $\mathcal{M}[\Tilde{\bf x}]$ and $\mathcal{G}[\Tilde{\bf x}]$,
\begin{equation}
    \mathcal{M}[\Tilde{\bf x}] = \int_{\tilde{\Gamma}} \tau\Tilde{\kappa}' \frac{\partial^2G(\Tilde{\bf x}-\Tilde{\bf x}')}{\partial {\bf n}(\Tilde{\bf x})\partial {\bf n}(\Tilde{\bf x}')}ds(\Tilde{\bf x}'), \label{MO}
\end{equation}
\begin{equation}
    \begin{split}
    \mathcal{G}[\Tilde{\bf x}] = \frac{1}{2}\frac{K_1+K_2}{K_1K_2} (1-\frac{\Tilde{A}}{\pi |\Tilde{\bf x}|^2})\Tilde{\bf x}\cdot{\bf n} + \frac{K_2-K_1}{K_1K_2} \int_{\tilde{\Gamma}} (1-\frac{\Tilde{A}}{\pi |\Tilde{\bf x}'|^2}) \Tilde{\bf x}'\cdot{\bf n} \frac{\partial G(\Tilde{\bf x} - \Tilde{\bf x}')}{\partial {\bf n}(\Tilde{\bf x})}ds(\Tilde{\bf x}') \\
    + \frac{K_2-K_1}{K_1K_2}\frac{\Tilde{A}}{\pi} \int_{\tilde{\Gamma}} \ln{|\Tilde{\bf x}'|} \frac{\partial^2G(\Tilde{\bf x}-\Tilde{\bf x}')}{\partial {\bf n}(\Tilde{\bf x})\partial {\bf n}(\Tilde{\bf x}')}ds(\Tilde{\bf x}'). \label{selfG}
    \end{split}
\end{equation}
These operators incorporate high-order derivatives of the interface. Moreover, at point $\Tilde{\bf x} = \Tilde{\bf x}'$, the kernel exhibits a logarithmic singularity. Evaluations of the singular integrals are explained in the next section.
Then, we rewrite Eq. (\ref{equaG}) as
\begin{equation}
\dot{R}R^2 = -\frac{\mathcal{M}[\Tilde{\bf x}]}{\mathcal{G}[\Tilde{\bf x}]} = C, \label{fluxc}
\end{equation}
where the time and space are separated, and $C$ is the flux constant. We rewrite it as
\begin{equation}
    \begin{split}
    \mathcal{M}[\Tilde{\bf x}] + C \mathcal{G}[\Tilde{\bf x}] = 0, \label{fluxc2}
    \end{split}
\end{equation}
which is a generalized nonlinear eigenvalue problem. The self-similar shape $\mathbf{\Tilde{x}}$ is the eigenvector and the flux constant $C$ is the eigenvalue. The operators $\mathcal{M}[\Tilde{\bf x}]$ and $\mathcal{G}[\Tilde{\bf x}]$ depend non-linearly and non-locally on the eigenvector ${\bf x}$. It is clear that circles satisfy Eq. (\ref{fluxc2}) with arbitrary number $C$. We are interested in noncircluar self-similar shapes.
Note that $\displaystyle\mathcal{M}[\beta \mathbf{x}]= \beta^{-2}\mathcal{M}[\mathbf{x}]$ and $\displaystyle\mathcal{G}[\beta \mathbf{x}]=\beta\mathcal{G}[\mathbf{x}]$ for any $\beta>0$..  If $\displaystyle(\mathbf{x},C)$ is an eigenpair, then $\displaystyle(\beta\mathbf{x},\frac{C}{\beta^3})$ is also an eigenpair. Thus, an eigenvector $\mathbf{x}$ can be paired with any given $C$ by adjusting its magnitude. We introduce a quasi-Newton method to solve Eq. (\ref{fluxc2}) with a pre-specified $C$ in the following section.

\section{Quasi-Newton scheme for self-similar shapes}
We use a quasi-Newton method to solve Eq. (\ref{fluxc2}). It is convenient to parameterize the interface by the polar angle $\alpha$ and solve for the Fourier modes in the polar angle representation. We specify the interface radius as a  combination of $cosine$ Fourier modes,
\begin{equation}
    \Tilde{r}(\alpha) = \sum_{k=0}^{N_1-1} \Tilde{\delta}_k \cos k\alpha,
\end{equation}
where $N_1$ is the total number of modes and $\Tilde{\delta}_k$ is the coefficient of the $k$th mode. Then the interface positions $(\Tilde{x}, \Tilde{y})$ are given as
\begin{eqnarray}
    \Tilde{x}(\alpha) = \Tilde{r}(\alpha)\cos \alpha, \quad
    \Tilde{y}(\alpha) = \Tilde{r}(\alpha)\sin \alpha.
\end{eqnarray}
The arc length variable $s$ in Eq. (\ref{fluxc2}) is related to $\alpha$ via $\displaystyle s_\alpha = \sqrt{\Tilde{x}_\alpha^2 + \Tilde{y}_\alpha^2}$. Note that $ds(\Tilde{x}) = s_\alpha d\alpha$.

Let ${\bf f}$ be a discretization of the left-hand side (LHS) of Eq. (\ref{fluxc2}). The discrete problem consists of finding $\Tilde{\delta}$'s for which
\begin{equation}
    {\bf f}(\Tilde{\delta}_0, \Tilde{\delta}_1, \Tilde{\delta}_2, \cdots, \Tilde{\delta}_{N_1-1}) = 0,
\end{equation}
at the interface node points $\alpha_i = i\Delta\alpha$ with $\Delta\alpha = 2\pi/N_2 $ for $i=0, \cdots, N_2-1$. Here $N_2$ is the total number of the interface node points.

It is clear that $\displaystyle \partial G(\Tilde{\bf x} - \Tilde{\bf x}')/\partial {\bf n}(\Tilde{\bf x}) = 1/(2\kappa(\Tilde{\bf x}))$ at $\Tilde{x}' = \Tilde{x}$. So the first integral of Eq. (\ref{selfG}) can be evaluated using a standard quadrature scheme. According to potential theory \cite{kress1989linear} and the Dirichlet-Neumann mapping \cite{greenbaum1993laplace}, the second integral of Eq. (\ref{selfG}) can be written as
\begin{equation}
    \int_{\Tilde{\Gamma}} \varphi(\Tilde{\bf x}') \frac{\partial^2G(\Tilde{\bf x}-\Tilde{\bf x}')}{\partial {\bf n}(\Tilde{\bf x})\partial {\bf n}(\Tilde{\bf x}')}ds(\Tilde{\bf x}') = \frac{d}{ds(\Tilde{\bf x})} \int_{\Tilde{\Gamma}} \frac{d\varphi}{ds}(\Tilde{\bf x}') G(\Tilde{\bf x}-\Tilde{\bf x}')ds(\Tilde{\bf x}') = \frac{1}{2\pi}\int_{\Tilde{\Gamma}} \frac{d\varphi}{ds}(\Tilde{\bf x}') \frac{(\Tilde{\bf x} - \Tilde{\bf x}')^\bot \cdot {\bf n}}{|\Tilde{\bf x} - \Tilde{\bf x}'|^2}ds (\Tilde{\bf x}'),
\end{equation}
where ${\bf x}^\bot = (y, -x)$. Thus, we use an alternating point trapezoidal rule to evaluate this integral and achieve spectral accuracy \cite{greenbaum1993laplace}. The integral in Eq. (\ref{MO}) can be similarly treated to handle the singularity. 

Our approach to address Eq. (\ref{fluxc2}) involves employing the well-known classical quasi-Newton method in conjunction with a line-search method \cite{press1992numerical}. Notably, the methodology bears striking similarities to that adopted by Li {\it et al.} \cite{li2004nonlinear} in his pursuit of determining the self-similar shape of a growing crystal.

The quasi-Newton method uses the iteration,
\begin{equation}
    \Tilde{\bm{\delta}}^{j+1} = \Tilde{\bm{\delta}}^{j} - {\bf J}^{-1} {\bf f}(\Tilde{\bm{\delta}}^{j}), \quad j=0,1,2,\cdots,
\end{equation}
where ${\bf J}$ is the Jacobian matrix, ${\bf J} = \nabla _{\Tilde{\delta}}{\bf f}$. Due to the highly nonlinear and nonlocal nature of  ${\bf f}$, establishing an explicit determination of  ${\bf J}$ is very difficult, as ${\bf f}$ relies heavily on the $\Tilde{\delta}_1, \Tilde{\delta}_2, \cdots, \Tilde{\delta}_{N_1 -1}$. Therefore, we use a finite difference approximation to calculate the matrix ${\bf J}$,
\begin{equation}
    J_{ij}=\frac{1}{h}\left(f_i(\Tilde{\delta}_0, \Tilde{\delta}_1, \Tilde{\delta}_2, \cdots, \Tilde{\delta}_j+h, \cdots,  \Tilde{\delta}_{N_1-1}) - f_i(\Tilde{\delta}_0, \Tilde{\delta}_1, \Tilde{\delta}_2, \cdots, \Tilde{\delta}_j, \cdots,  \Tilde{\delta}_{N_1-1})\right)
\end{equation}
where $h$ is small. The quasi-Newton iteration is performed until $\max|{\bf f}|$ is less than a prescribed tolerance. Once $\mathbf{\Tilde{x}}$ is obtained, we scale the vector with $\delta_0=1$ and calculate the corresponding flux constant $\displaystyle C=-\frac{\mathcal{M}[\mathbf{\Tilde{x}}]}{\mathcal{G}[\mathbf{\Tilde{x}}]}$.

\section{Self-Similar shapes}
We implement the quasi-Newton scheme to solve Eq.(\ref{fluxc2}) with $\tau=1$. We calculate a shape factor, 
\begin{equation}
    \frac{\delta}{\bf R} = \max||\Tilde{\bf x}|/\Tilde{\bf R}_\mathrm{eff}-1|,
\end{equation}
where $\Tilde{\bf x}$ is the position vector of the interface and  $\Tilde{\bf R}_\mathrm{eff}$ is the effective radius of the nonlinear (self-similar) shape. For the initial configuration of the quasi-Newton scheme, $\Tilde{\delta}_0=1$ is set. 

\subsection{Reliability of method}

First, we give a resolution study of a 3-fold dominant self-similar shape to assess the accuracy of the quasi-Newton solver when $N_1=128$. We use four resolutions, $N_2=256$, $512$, $1024$, and $2048$ with the initial $C_0=30$ and the initial guess $\displaystyle\Tilde{\delta}_3 = 0.2$ ( other modes are zeros except for $\Tilde{\delta}_0 = 1$ ) to compute the self-similar shapes. All cases experience a 3-fold self-similar shape. A Fourier analysis shows that only mode 3 and its harmonics present nonzero Fourier coefficients. 
Our results show that $\displaystyle \frac{\delta}{\bf R} = 0.243225$, $0.243224$, $0.243223$, and $0.243219$, for $N_2=256$, $512$, $1024$, and $2048$, respectively. 
Thus, we obtain the quasi-Newton scheme to be spectrally accurate with the following error expansion,
\begin{equation}
     \frac{\delta}{\bf R} = \left( \frac{\delta}{\bf R} \right)^* + \beta_1 \mathrm{e}^{-\beta_2N_2},
\end{equation}
where $\displaystyle \left( \frac{\delta}{\bf R} \right)^*= 0.2431865$ is the exact solution, $\beta_1 = 3.961\times 10^{-5}$, and $\beta_2 = 9.464\times 10^{-5}$. The error is the difference between $\displaystyle \frac{\delta}{\bf R}$ and exact solution $\displaystyle \left( \frac{\delta}{\bf R} \right)^*$. 
\begin{figure}[h]
    \centering    
    \includegraphics[width=0.45\linewidth]{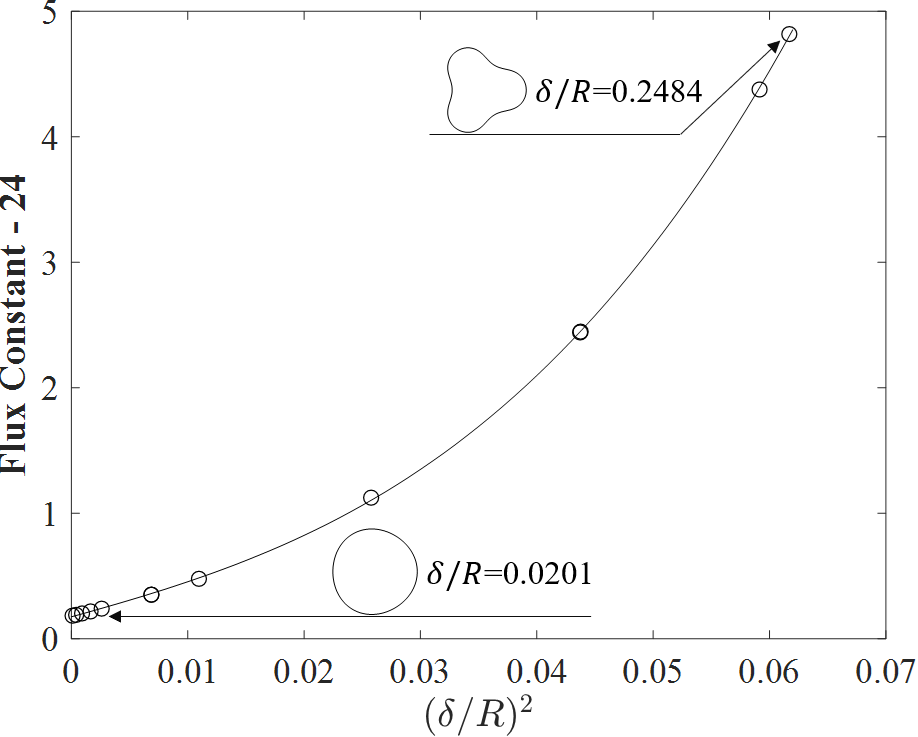}(a)
    \includegraphics[width=0.45\linewidth]{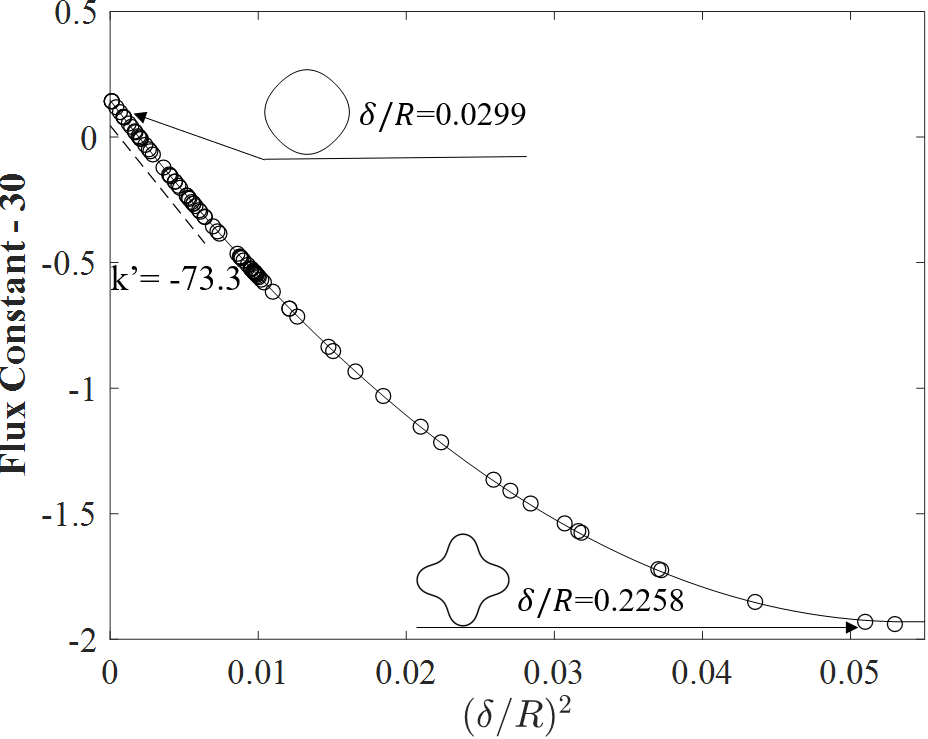}(b)
    \caption{The difference of flux constants between the linear theory and nonlinear theory for (a) 3-fold dominant self-similar shapes and (b) 4-fold dominant self-similar shapes.} 
    \label{fold34}
\end{figure}
Then, we study the difference of flux constants between the linear theory and nonlinear theory for 3-fold and 4-fold dominant self-similar shapes. 
Figures \ref{fold34}(a) and \ref{fold34}(b) show the different nonlinear results that obtained by varying $C_0$, $\Tilde{\delta}_3$ and $\Tilde{\delta}_4$ for the initial guess. The nonlinear flux constants are obtained form Eq. (\ref{fluxc}). The linear flux constant is $24$ and $30$ by $\displaystyle C = \frac{k(k^2-1)}{k-2}$ with $k=3$ and $k=4$, respectively. The illustrations in Fig. \ref{fold34} also show the associated interface shapes for different values of $\delta/R$. For the 3-fold dominant self-similar shape, the deviation of the nonlinear results from linear theory is cubical in $(\delta/R)^2$. Nonlinear effects result in an elevated flux constant. For the 4-fold dominant self-similar shape, the deviation of the nonlinear results from linear theory is quadratic in $\delta/R$, when $\delta/R$ is small ($\delta/R \leq 0.08$). In other words, the deviation is linear in $(\delta/R)^2$, with a slope of $\text{k'}=-73.3$. The deviation of the nonlinear results from linear theory is quadratic in $(\delta/R)^2$ for large $\delta/R$. Nonlinear effects result in a lowered flux constant for the 4-fold dominant self-similar shape. Analogous results with a decay behaviour are found for other k-fold ($k>4$) self-similar shapes as well.

\begin{figure}[t]
    \centering    
    \includegraphics[width=0.75\linewidth]{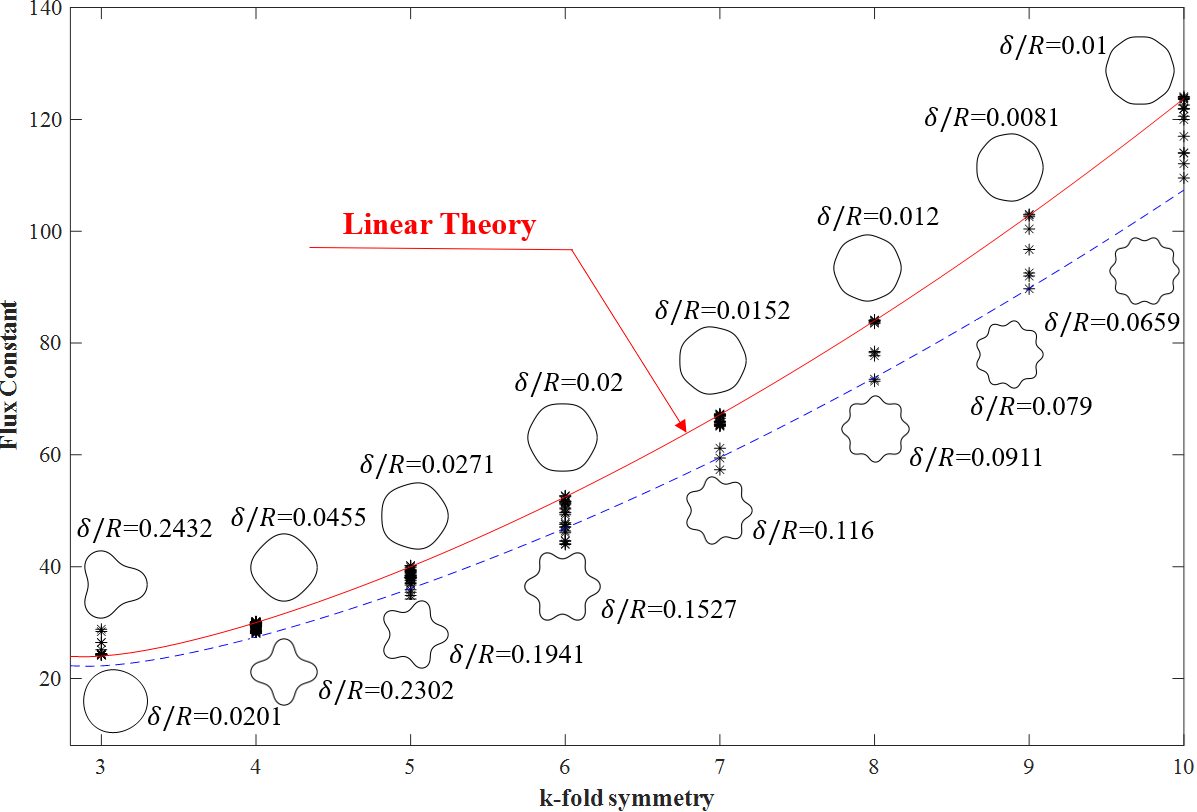}
    \caption{Flux constants of self-similar shapes and selected morphologies. Flux constants from linear theory (solid line) are given by $\displaystyle C = \frac{k(k^2-1)}{k-2}$. Symbols denote the nonlinear results. The dashed line represents the best fit for the nonlinear self-similar shapes with a large shape factor $\delta/R$, $\displaystyle C = \frac{k(k^{1.939}-1)}{k-2}$ ($k\ge 4$).} 
    \label{kfold}
\end{figure}

We consider the nonlinear effects for different symmetries. Figure \ref{kfold} shows the flux constants and morphologies for general $k$-fold self-similar shapes. We plot the linear theory prediction $\displaystyle C = \frac{k(k^2-1)}{k-2}$ as the solid curve (red), while symbols denote the nonlinear results. The flux constants match well with the linear theory predictions when the perturbations are small. The dashed line (blue) represents the best fit for the nonlinear self-similar shapes with a large shape factor $\delta/R$, $\displaystyle C = \frac{k(k^{1.939}-1)}{k-2}$ with $k\ge 4$. The flux constants are smaller than the results of the linear theory when the perturbations are big, while the 3-fold case behaves conversely.

\subsection{The influence of initial parameters}

\subsubsection{The initial guess}
\begin{figure}[h]
    \centering    
    \includegraphics[width=0.45\linewidth]{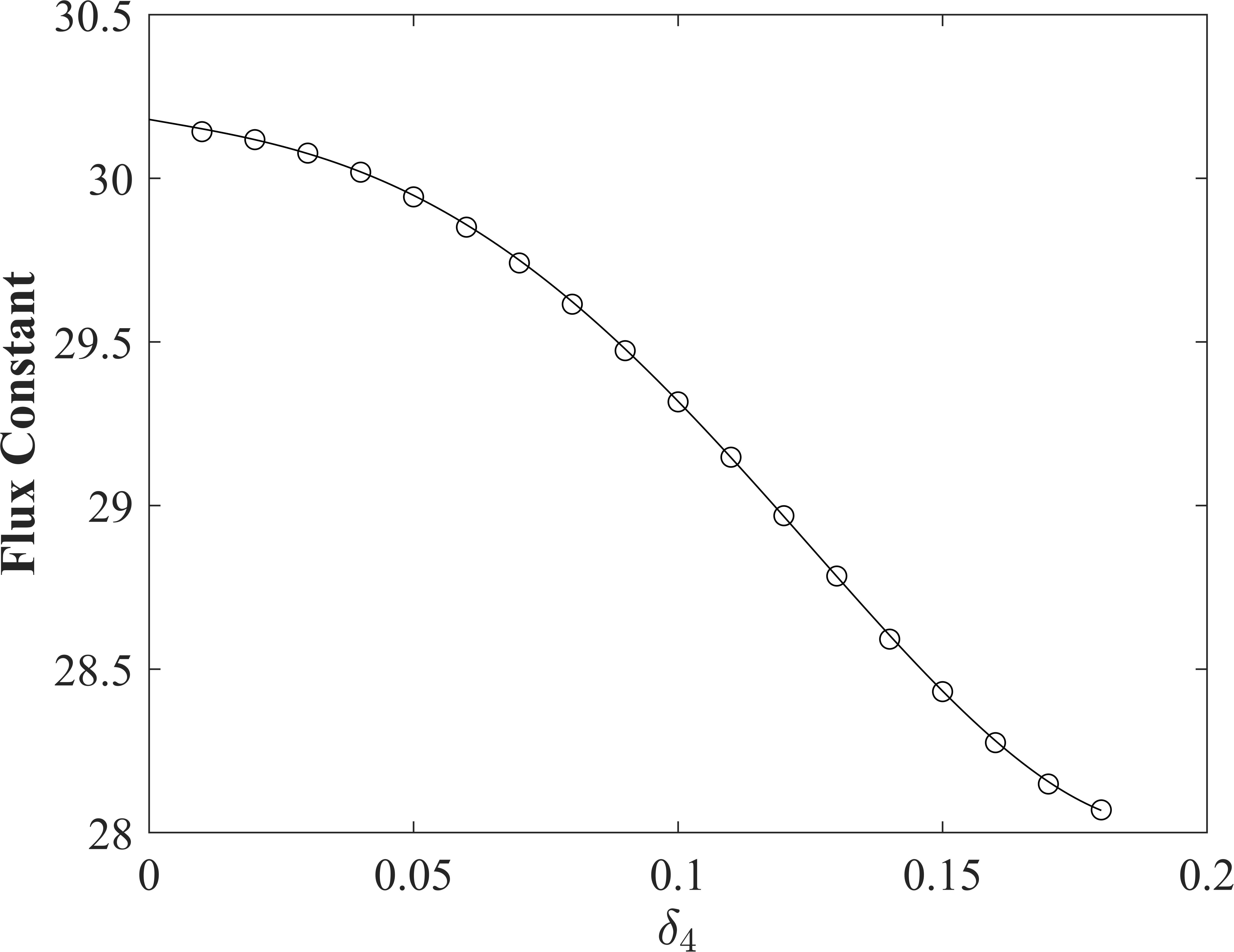}(a)
    \includegraphics[width=0.45\linewidth]{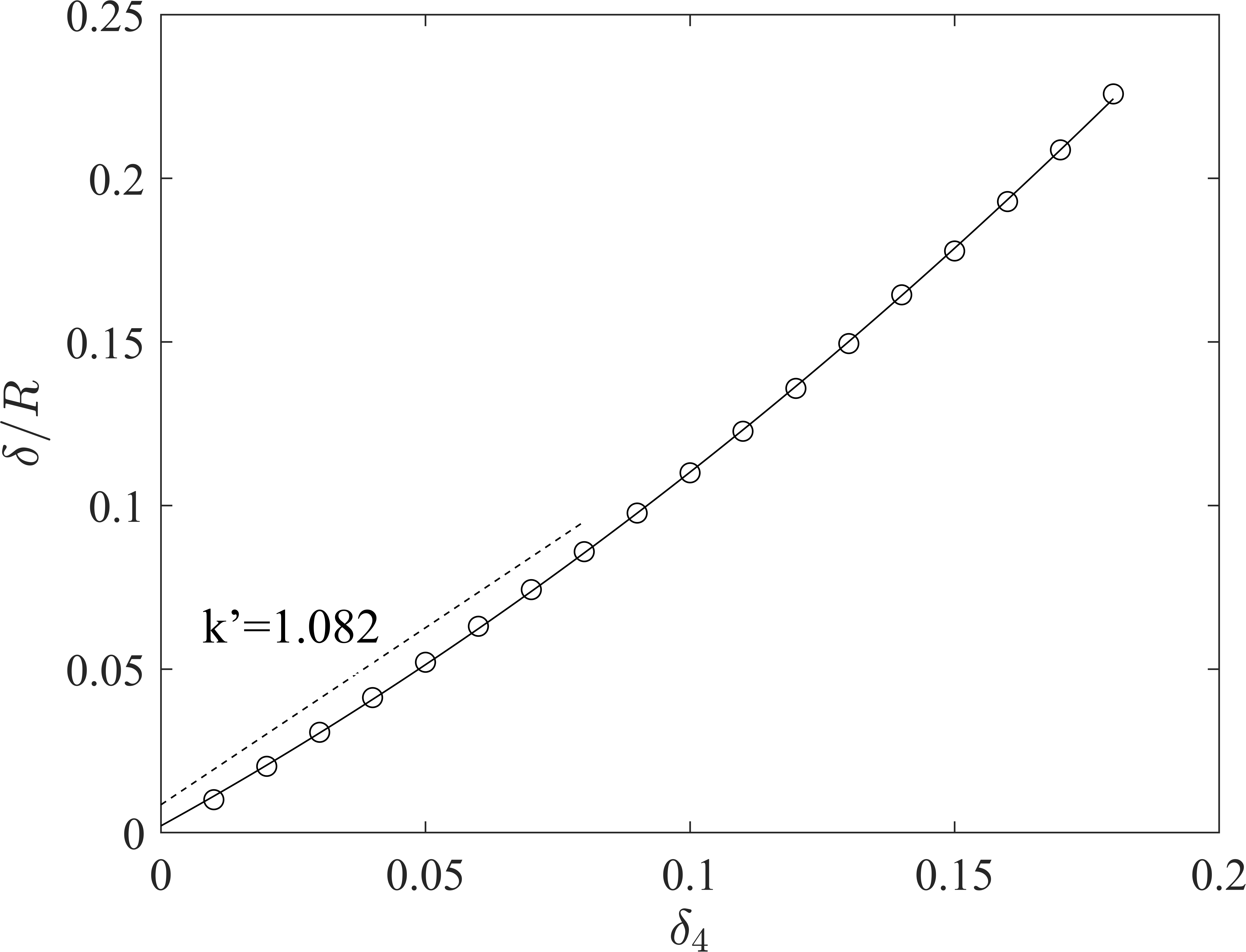}(b)
    \caption{The effect of the initial guess of a single mode on (a) the flux constant and (b) the shape factor $\delta/R$. We set $C_0 = 30$. The different nonlinear simulations are obtained by varying $\Tilde{\delta}_4$ in the initial guess.} 
    \label{IniSingGue}
\end{figure}

We discuss the effect of initial guesses on self-similar shapes. First, we examine the case where the initial guess consists of a single mode, e.g. $\Tilde{\delta}_0=1$, $\Tilde{\delta}_4\neq 0$, and all other modes perturbation is zero. We set $C_0 = 30$ and change the initial value of $\Tilde{\delta}_4$ to get the evolution of the flux constant and the shape factor. In Fig. \ref{IniSingGue}(a), the flux constant is quadratic in the initial guess $\Tilde{\delta}_4$, when $\Tilde{\delta}_4$ is small ($\Tilde{\delta}_4 \leq 0.08$). It is quartic in $\Tilde{\delta}_4$ for large initial guesses. In Fig. \ref{IniSingGue}(b), the shape factor $\delta/R$ is linear in the initial guess $\Tilde{\delta}_4$ with a slope of $\text{k'}=1.082$, when $\Tilde{\delta}_4$ is small. It is quadratic in $\Tilde{\delta}_4$ for large initial guess. These results are consistent with the relationship between the flux constant and the shape factor in the previous section.

\begin{figure}[h]
    \centering    
    \includegraphics[width=0.45\linewidth]{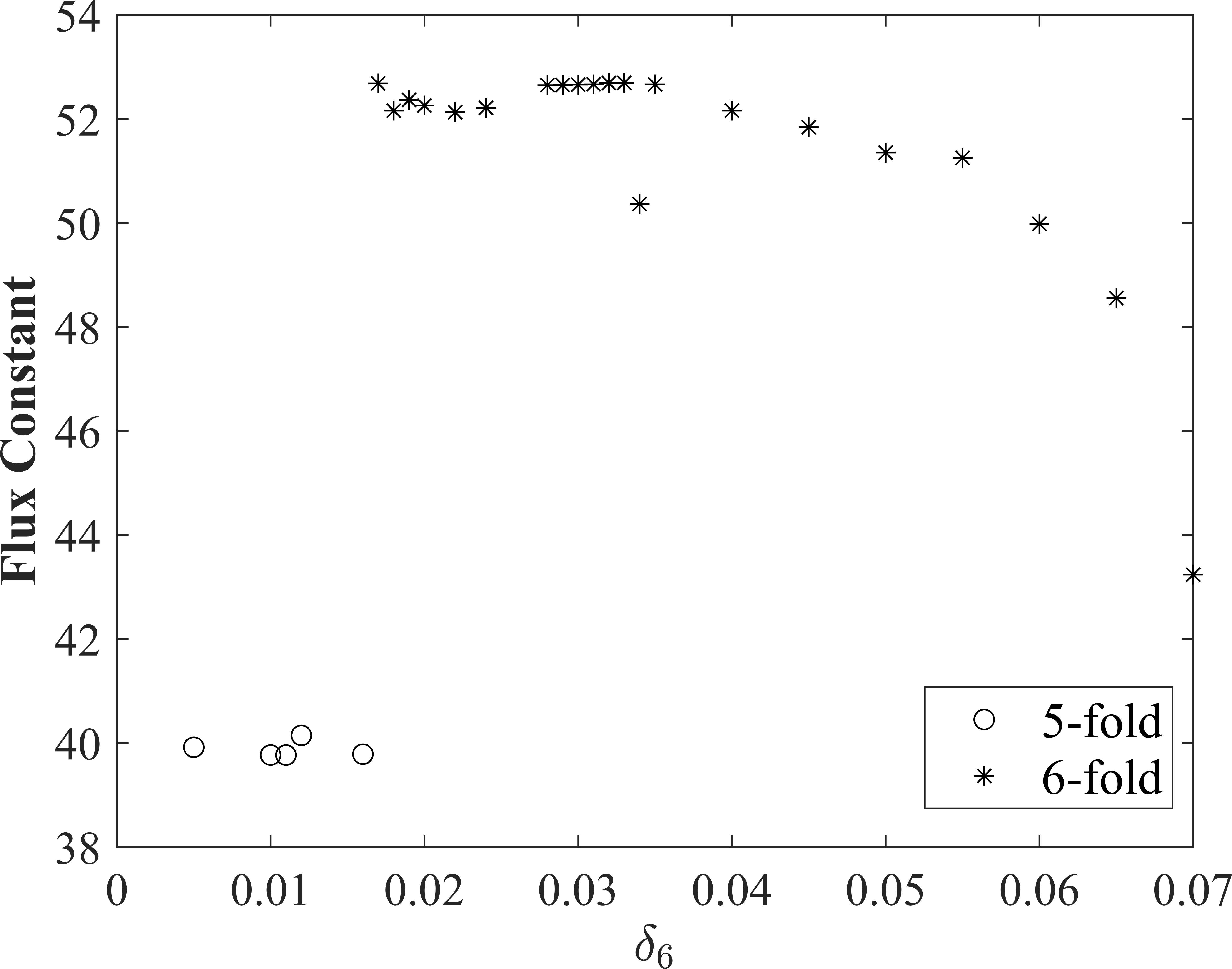}
    \caption{The effect of the initial guess with mixed modes on the flux constant. We set $C_0 = 50$ and $\Tilde{\delta}_5 = 0.1$. The different nonlinear simulations are obtained by varying $\Tilde{\delta}_6$.
    } 
    \label{IniMixGue}
\end{figure}

Then, we examine the case where the initial guess includes mixed modes, e.g. $\Tilde{\delta}_0=1$, $\Tilde{\delta}_5\neq 0$, and $\Tilde{\delta}_6\neq 0$. We set $C_0 = 50$ and $\Tilde{\delta}_5 = 0.1$. The evolution of the flux constant with the change of $\Tilde{\delta}_6$ is shown in Fig. \ref{IniMixGue}. When $\Tilde{\delta}_6 \leq 0.016$, we obtain a 5-fold dominant self-similar shape. When $\Tilde{\delta}_6 > 0.016$, we obtain a 6-fold dominant self-similar shape. 
The flux constant gradually decreases as the initial guess $\Tilde{\delta}_6$ gradually increases. That is to say, the nonlinearity of the self-similar shape gradually increases.

\subsubsection{The initial parameter $C_0$}

\begin{figure}[h]
    \centering    
    \includegraphics[width=0.9\linewidth]{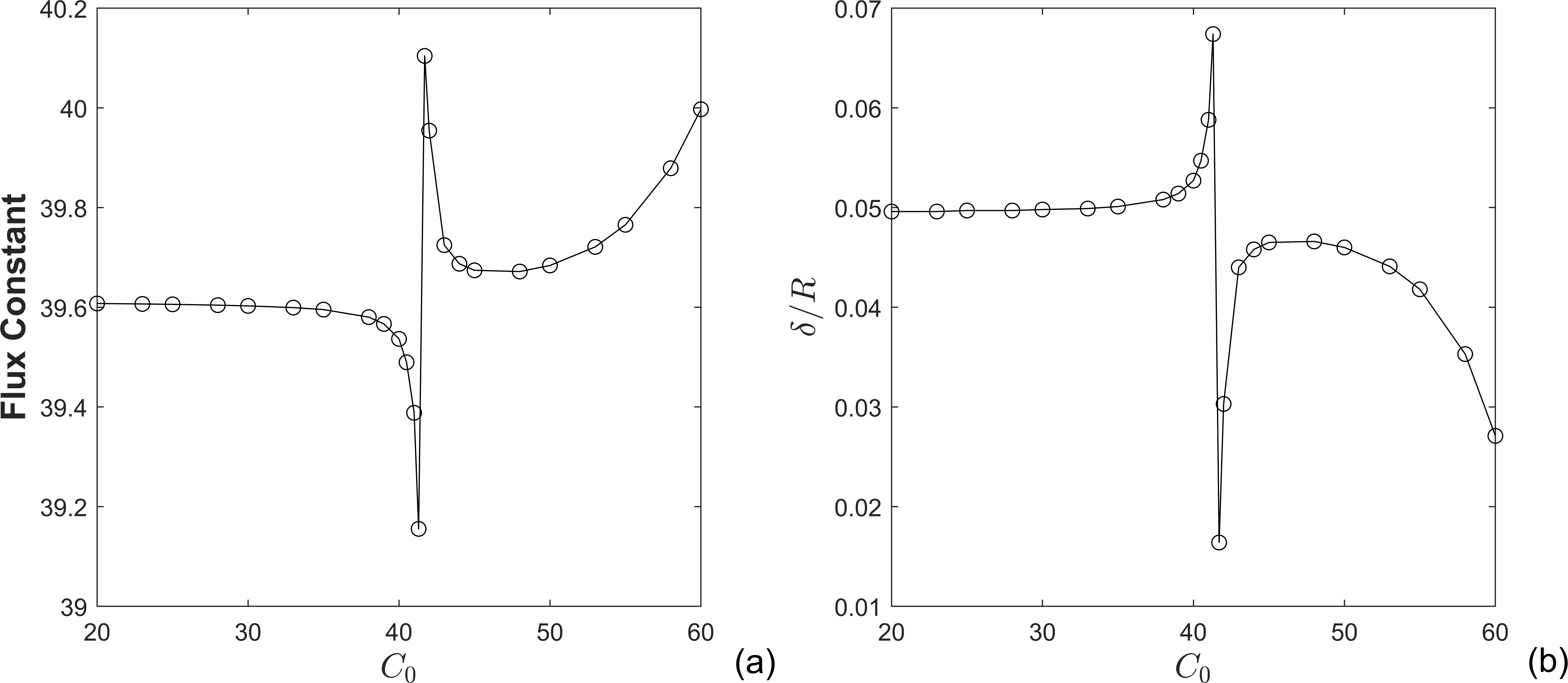}
    \caption{The effect of $C_0$ with initial guesses of a single mode. We set the initial guess $\Tilde{\delta}_5 = 0.05$. The different nonlinear simulations are obtained by varying $C_0$. } 
    \label{CSingGueCsh}
\end{figure}

In this subsection, we study the effect of the initial parameter $C_0$ on self-similar shapes. We set the initial guess to be a single mode $\Tilde{\delta}_5 = 0.05$. We vary $C_0$ to get the evolution of the flux constant and the shape factor. The results are shown in Figs. \ref{CSingGueCsh}(a) and \ref{CSingGueCsh}(b), respectively. As $C_0$ increases, the flux constant changes little at first, and starts to decrease rapidly at about $C_0=35$. When the flux constant is small to a certain extent $C^*$, there is a jump increment. The flux constant becomes large, and then decreases. Later, it increases gradually. 
Such a jump may be due to the presence of nonunique solutions in Eq. (\ref{fluxc2}). Roughly speaking, the LHS of Eq. (\ref{fluxc2}) may achieve its local extreme at this initial configuration with $C^*$. 
 When $C_0 \leq C^*$, one solution is obtained and when $C_0 > C^*$, another solution is obtained. These two solutions exhibit significant differences.

\begin{figure}[h]
    \centering    
    \includegraphics[width=0.45\linewidth]{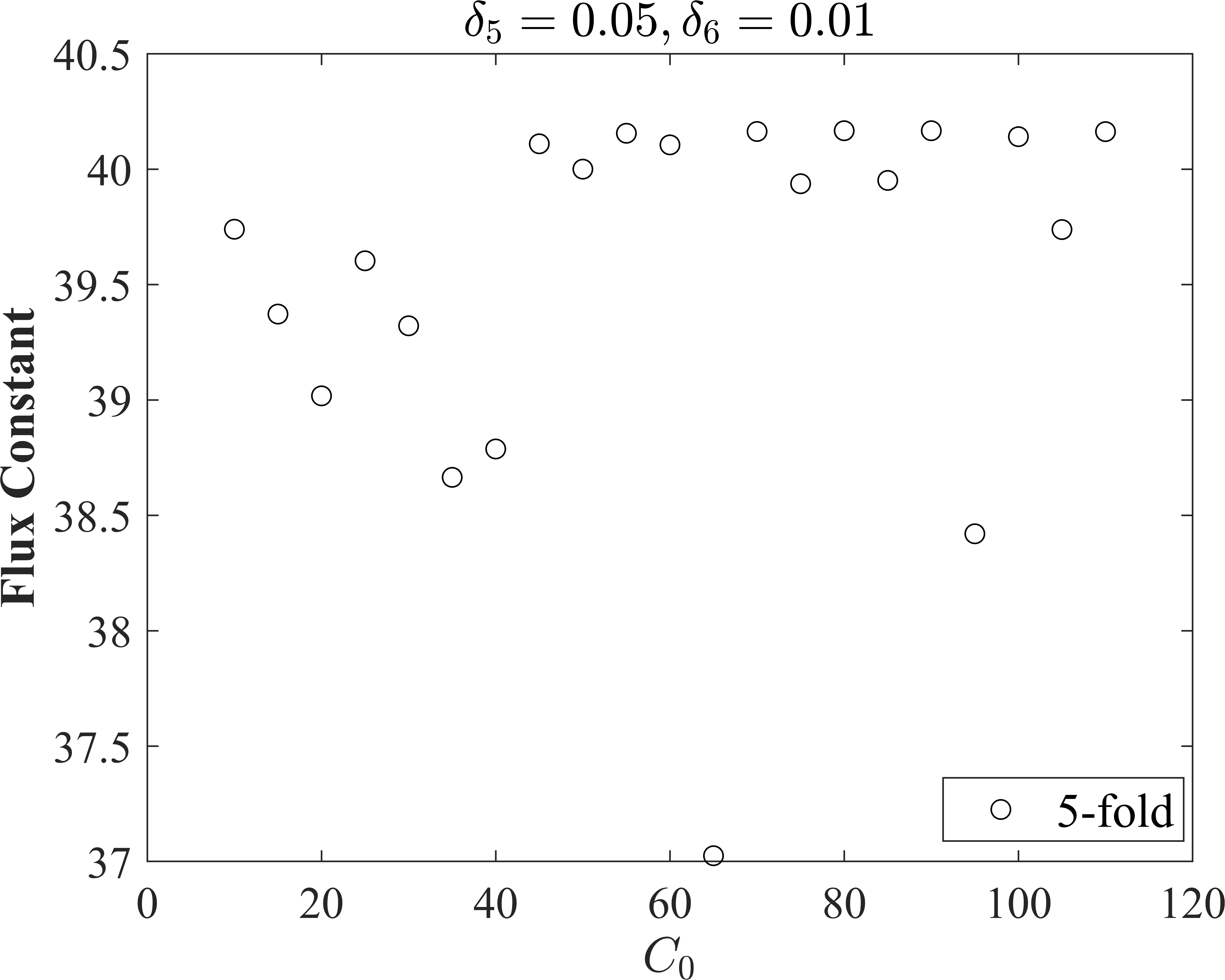}(a)
    \includegraphics[width=0.45\linewidth]{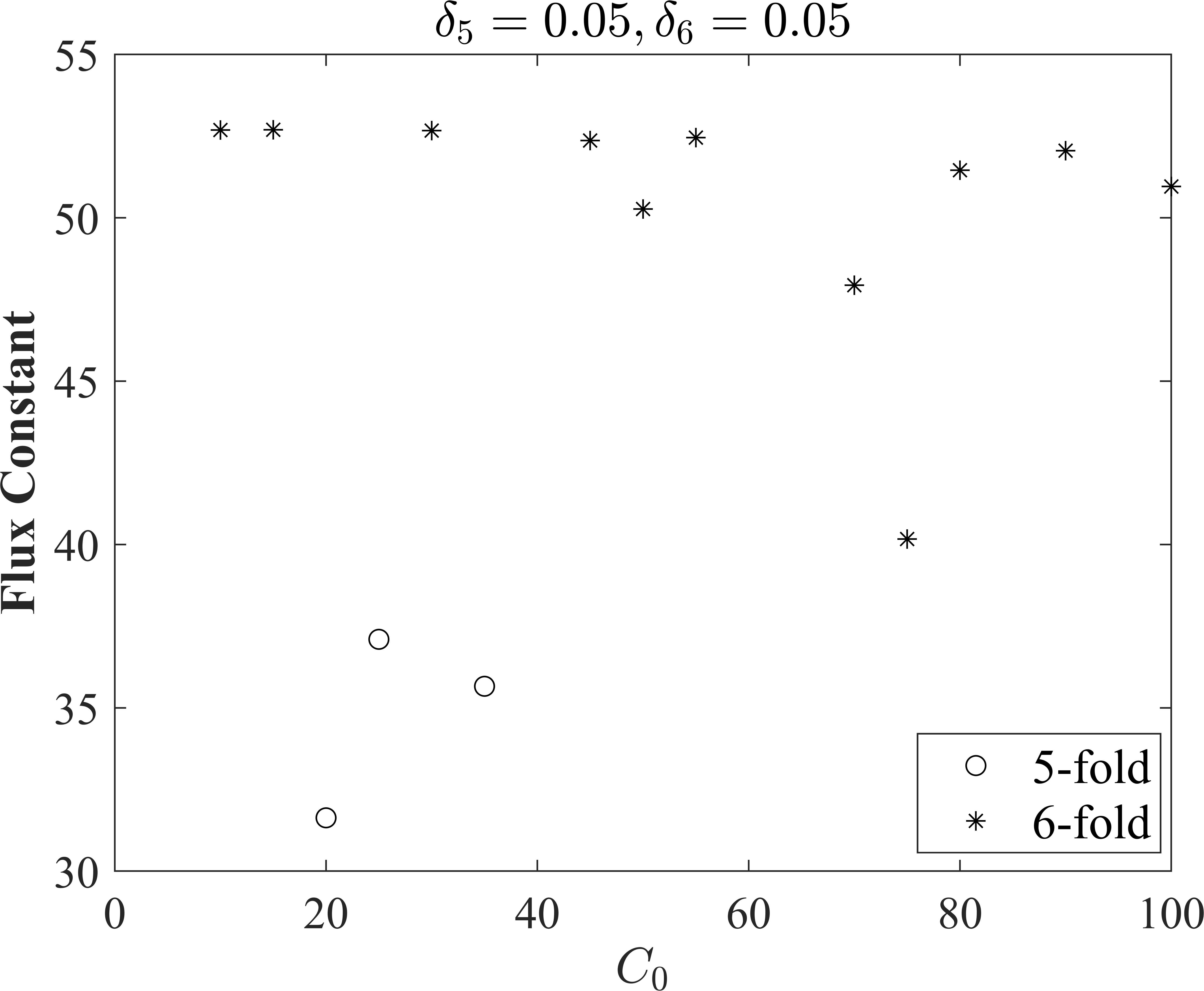}(b)
    \caption{The effect of $C_0$ with initial guesses of mixed modes on the flux constant. The different nonlinear simulations are obtained by varying $C_0$. (a) We set the initial guesses $\Tilde{\delta}_5 = 0.05$ and $\Tilde{\delta}_6 = 0.01$. (b) We set the initial guesses $\Tilde{\delta}_5 = 0.05$ and $\Tilde{\delta}_6 = 0.05$.} 
    \label{CMixGue}
\end{figure}

Moreover, we investigate the influence of $C_0$ on self-similar shapes in two cases with initial guesses $\Tilde{\delta}_5\neq 0$ and $\Tilde{\delta}_6\neq 0$. In one case, the initial perturbation of these two modes differs significantly ( $\Tilde{\delta}_5 = 0.05$ and $\Tilde{\delta}_6 = 0.01$ are set as shown in Fig. \ref{CMixGue}(a)). It shows that the flux constant fluctuates with the increase of the initial $C_0$. The range of the change in flux constant is small. All we get are 5-fold dominant self-similar shapes. When one perturbation is much larger than the other, the self-similar shape is not highly affected by $C_0$ and always dominated by the larger one. In the other case,  the initial perturbations of the two modes are equal ($\Tilde{\delta}_5 = \Tilde{\delta}_6 = 0.05$ are set as shown in Fig. \ref{CMixGue}(b)). It shows a wide range of changes in the flux constant as the initial $C_0$ increases. Self-similar shapes also bounce back and forth between five folds and six folds. In other words, when the values of the two modes in the initial guesses are equal, the self-similar shape is greatly affected by the initial $C_0$.

\subsection{Non-trivial shapes}
Our formulation and scheme are nonlinear and non-trivial. Trivially, given an initial guess $\Tilde{\delta}_k\neq 0$, only a k-fold dominant self-similar shape can be obtained. However, our formulation allows for the computation of self-similar shapes with harmonics of $k$. For example, when the initial guess is set as $\Tilde{\delta}_5=0.3$ and $\Tilde{\delta}_8=0.1$, the scheme exhibits 10-fold and 16-fold dominant self-similar shapes, respectively. 

\begin{figure}[h]
    \centering    
    \includegraphics[width=0.7\linewidth]{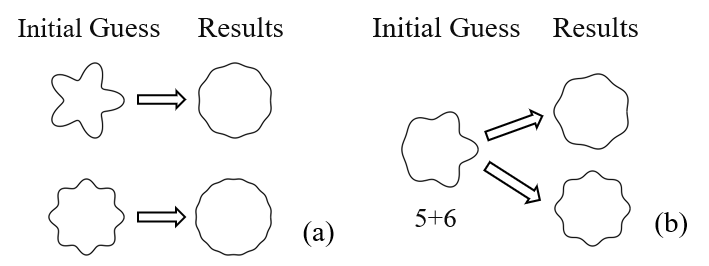}
    \caption{Diagram of non-trivial computing cases: (a) computation of self-similar shapes of harmonic order; (b) computation of self-similar shapes not included in the initial guesses.} 
    \label{mutotal}
\end{figure}

Furthermore, when the initial guess contains mixed modes, we are able to obtain self-similar shapes that are not initially included. For example, we set $C_0=65$, $\Tilde{\delta}_5 =\Tilde{\delta}_6=0.05$ and obtain a 7-fold dominant self-similar shape. When we modify the initial configuration as $C_0=60$ and $\Tilde{\delta}_5 =\Tilde{\delta}_6=0.1$, we get a 8-fold dominant self-similar shape instead. In Fig. \ref{mutotal}, we illustrate these non-trivial computation cases. When we choose the initial guess and initial $C_0$ properly, we may calculate any self-similar shape. This goes beyond predictions of linear theory and experiences similar results with fully nonlinear simulations \cite{PhysRevLett.102.174501, zhao2017efficient}.


\section{Conclusion}
In summary, we have developed a nonlinear theory for the self-similar interface between two immiscible fluids in a radial Hele-Shaw cell. A generalized nonlinear eigenvalue problem  $\mathcal{M}[{\bf x}] + C \mathcal{G}[{\bf x}] = 0$ is obtained. Here the self-similar shape ${\bf x}$ is the eigenvector and the flux constant $C$ is the eigenvalue. The problem is challenging to solve due to its highly nonlocal, nonlinear nature, and the presence of singularities. 
We have investigated a quasi-Newton scheme to calculate the self-similar solution. Our results indicate that nonlinear, noncircular self-similarly interfaces indeed exist. Nonlinearity reduces the flux constant of k-fold dominant self-similar shapes ($k\ge4$) compared to linear theory predictions. However, the flux constant of the 3-fold dominant self-similar shape surpasses the linear theory predictions. 

Both the initial guess of $\Tilde{\delta}_k$ and the initial parameter $C_0$ have an impact on the self-similar shapes. 
As the initial guess of $\Tilde{\delta}_k$ increases, nonlinearity of the system is enhanced and the self-similar solutions deviate from linear predictions.
When the initial guess contains mixed modes, different initial guesses may yield self-similar shapes of different folds. The influence of the initial parameter $C_0$ on the self-similar shape is more complex. In particular, as long as the initial guesses are appropriately controlled, any fold of self-similar shape may be carried out. For example, we are able to obtain both a 7-fold and 8-fold dominant self-similar shape from the initial guesses with mode 5 and mode 6 by changing $C_0$. This behaviour transcends the boundaries of linear theory and attains results on par with simulations. The nonlinear theory bridges the linear theory  and simulations.

This paper is primarily focused on fluid dynamics in a Hele-Shaw cell. However, the research framework established here can be easily extended to include gravitational or electrical fields. This expansion opens up opportunities for future investigations, which we intend to pursue. 





\section*{CRediT authorship contribution statement}
\textbf{Wang Xiao:} Conceptualization, Methodology, Software, Formal analysis, Investigation, Writing - Original Draft. \textbf{Lingyu Feng:} Writing - Review \& Editing, Visualization. \textbf{Kai Liu:} Resources, Writing - Review \& Editing, Supervision, Funding acquisition. \textbf{Meng Zhao:} Conceptualization, Methodology, Validation, Formal analysis, Resources, Writing - Review \& Editing, Supervision, Funding acquisition.

\section*{Declaration of competing interest}
The authors declare that they have no known competing financial interests or personal relationships that could have appeared to influence the work reported in this paper.

\section*{Data availability}
Data will be made available on request.

\section*{Acknowledgments}
This research was supported by the National Natural Science Foundation of China 12301553. M. Z. also thanks the support from the Ministry of Education Key Lab in Scientific and Engineering Computing.

\bibliography{cite}

\end{document}